\documentclass[11pt, a4paper]{article}
\usepackage{amsmath,amssymb,latexsym,graphics,epsfig,hyperref}

\newtheorem{theo}{Theorem}[section]

\newtheorem{proposicio}[theo]{Proposition}

\def\prova{{\boldmath  $Proof.$}\hskip 0.3truecm}

\def\final{\mbox{ \quad $\Box$}}

\def\A{\mbox{\boldmath $A$}}

\def\N{\mathbb N}
\def\Re{\mathbb R}

\def\dist{\mathop{\rm dist }\nolimits}
\def\dgr{\mathop{\rm dgr }\nolimits}
\def\excu{\mbox{$\varepsilon_u$}}
\def\exc{\mbox{$\varepsilon$}}

\def\Ker{\mathop{\rm Ker }\nolimits}

\def\e{\mbox{\boldmath $e$}}
\def\j{\mbox{\boldmath $j$}}

\def\vecnu{\mbox{\boldmath $\nu$}}
\def\vecrho{\mbox{\boldmath $\rho$}}

\def\vecalpha{\mbox{\boldmath $\alpha$}}
\def\vece{\mbox{\boldmath $e$}}

\def\A{\mbox{\boldmath $A$}}

\def\E{\mbox{\boldmath $E$}}

\def\G{\Gamma}

\def\I{\mbox{\boldmath $I$}}

\def\J{\mbox{\boldmath $J$}}
\def\M{\mbox{\boldmath $M$}}
\def\N{\mbox{\boldmath $N$}}

\def\S{\mbox{\boldmath $S$}}

\def\ev{\mbox{\rm ev}}
\def\sp{\mbox{\rm sp}}

\def\Npetita{\mbox{{\tiny $N$}}}

\def\subN{\Npetita}

\def\dgr{\mathop{\rm dgr }\nolimits}
\def\ecc{\mathop{\rm ecc }\nolimits}

\def\Ker{\mathop{\rm ker }\nolimits}
\def\tr{\mathop{\rm tr }\nolimits}

\def\ev{\mathop{\rm ev }\nolimits}
\def\sp{\mathop{\rm sp }\nolimits}

\begin{document}
\title{On Some Approaches to the Spectral\\ Excess Theorem for Nonregular Graphs.
% Related Results and Open Problems
%\thanks{Research supported by the Ministerio de Educaci\'on y
%Ciencia, Spain, and the European Regional Development Fund under
%project MTM2008-06620-C03-01 and by the Catalan Research Council
%under project 2009SGR1387.}
}

\author{M.A. Fiol
\\ \\
{\small Universitat Polit\`ecnica de Catalunya}\\
%{\small Tilburg University} \\
{\small Departament de Matem\`atica Aplicada IV} \\
{\small Barcelona, Catalonia} \\
{\small e-mail:  {\tt fiol@ma4.upc.edu}}
}

%\date{}

\maketitle

\begin{abstract}
The Spectral Excess Theorem ({\em SPET\/}) for distance-regular graphs states that a regular (connected) graph is distance-regular if and only if its spectral-excess equals its average excess.  Recently, some local or global approaches to the {\em SPET\/} have been used to obtain new versions of the theorem for nonregular graphs, and also to study the problem of characterizing the graphs which have the corresponding distance-regularity property.
In this paper, some of these versions are related and compared, and some of their results are improved. As a result, a sufficient condition for a graph to be distance-polynomial is obtained.
\end{abstract}

\vskip 2cm

\section{Introduction}

The spectral excess theorem \cite{fg97} ({\em SPET\/}) states that a regular (connected) graph $\G$ is distance-regular if and only if its spectral-excess (a number which can be computed from the spectrum of $\G$) equals its average excess (the mean of the numbers of vertices at maximum distance from every vertex), see \cite{vd08,fgg10} for short proofs.
Since the paper \cite{fg97} appeared, different approaches (local or global) of the {\em SPET} have been given.
The interest of the  inequalities so obtained is the characterization of some kind of distance-regularity, and this happens when equalities are attained.

One of such versions concerns with the so-called pseudo-distance-regularity \cite{fgy96}, which is a natural generalization, for nonregular graphs, of the standard distance-regularity \cite{b93,bcn89}.
As shown in \cite{fgy96}, the price to be paid for speaking about distance-regularity in nonregular graphs is locality. More precisely, what is called ``pseudo-distance-regularity" around a vertex. Thus, generalizing a result of Godsil and Shawe-Taylor \cite{gst87}, the author  recently proved that  when the pseudo-distance-regularity is shared by all the vertices, the graph is either distance-regular or distance-biregular, see \cite{f12}.

In this paper the commonly used approaches to the {\em SPET\/}, that is, local and global, are compared and related. In particular,
Lee and Weng \cite{lw11} recently derived an inequality for nonregular graphs which is similar to the one that leads to the {\em SPET}, and posed the problem of characterizing the graphs for which the equality is attained. As a result of our study, we show that, in some cases,  such an equality is a sufficient condition a the graph $\G$ to be distance-polynomial, which requires that every distance matrix of $\G$ is a polynomial in its adjacency matrix (see Weichel \cite{w82}).

%%%%%%%%%%%%%%%%%%%%%%%%%%%%%%%%%%%%%%%%%%%%%%%%
\section{Preliminaries}
%%%%%%%%%%%%%%%%%%%%%%%%%%%%%%%%%%%%%%%%%%%%%%%%

\subsection{Graphs and their spectra}

Let us first recall some basic notation and results on which our study is based. For more background on spectra of graphs,
distance-regular graphs, and their characterizations, see
\cite{b93,bcn89,bh12,cds82,dkt12,f02,g93}.

Throughout this paper, $\G$ denotes a (finite, simple and connected) graph
with vertex set $V$, order $n=|V|$, and adjacency matrix
$\A$. The {\it distance} between two vertices $u$ and $v$ is
denoted by $\dist (u,v)$, so that the {\it eccentricity} of
vertex $u$ is $\exc_u=\max_{v\in V}\dist (u,v)$ and the {\it
diameter} of the graph is $D=\max_{u\in V}\excu$. The set of
vertices at distance $i$, from a given vertex $u\in V$ is
denoted by $\G_i(u)$, for $i=0,1,\dots,D$, and $N_i(u)=\G_0(u)\cup\G_1(u)\cup\cdots\cup \G_i(u)$. The degree of
vertex $u$ is $\delta_u=|\G(u)|=|\Gamma_1(u)|$.

The spectrum of $\G$ is denoted by
$
\sp \G = \sp \A = \{\lambda_0^{m_0},\lambda_1^{m_1},\dots,
\lambda_d^{m_d}\},
$
where the different eigenvalues of $\G$ are in decreasing order,
$\lambda_0>\lambda_1>\cdots >\lambda_d$, and the superscripts
stand for their multiplicities $m_i=m(\lambda_i)$.
In
particular, note that  $m_0=1$ (since $\G$ is
connected) and $m_0+m_1+\cdots+m_d=n$. Moreover, $\lambda_0$ has a positive eigenvector (usually called the {\em Perron vector\/}),
denoted by $\vecnu$ when it is normalized in such a way that its minimum component is
$1$ \cite{fgy96}, and $\vecalpha$ when the normalization condition is $\|\vecalpha\|^2=n$ \cite{lw11}.
For instance, if $\G$ is regular, we have $\vecnu=\vecalpha=\j$, the all-$1$ vector.

Let $\E_i$, $i=0,1,\ldots,d$, be the idempotents of $\A$, that is $\E_i=\frac{1}{\phi_i}\prod_{j\neq i}(\A-\lambda_j\I)$, where $\phi_i=\prod_{j\neq i}(\lambda_i-\lambda_j)$.
The {\em $u$-local multiplicity} of $\lambda_i$ is $m_{u}(\lambda_i)=\|\E_i\vece_{u}\|^2=(\E_i)_{uu}$ and  $\ev_{u}\G$ denotes the set of {\em $u$-local eigenvalues}, that is, those eigenvalues of $\G$ with nonzero $u$-local multiplicity. Since $\vece_{u}$ is unitary, we have $\sum_{i=0}^d m_{u}(\lambda_i)=1$. We refer to the pair $(\ev_u\G, m_{u})$, constituted by the set of $u$-local eigenvalues and the normalized weight function $m_{u}$ defined by the $u$-local multiplicities, as the {\em $u$-local spectrum} of $\G$.
Moreover, as $\E_0\vece_{u}=\frac{1}{\|\vecalpha\|^2}\langle\vece_{u}, \vecalpha \rangle \vecalpha
=\frac{\alpha_u}{n}\vecalpha$, we have that  $m_{u}(\lambda_0)=\frac{\alpha_u^2}{n}$ and, hence, $\lambda_0\in \ev_{u}\G$. Let $\ev^{\star}_{u}\G=\ev_u \G\setminus\{\lambda_0\}$ and $d_{u}=|\ev^{\star}_{u}\G|$. In \cite{fg99} it was shown that then the eccentricity of $u$ satisfies $\excu\le d_{u}$. When equality is attained, we say that $u$ is an {\em extremal vertex}.

Recall that, for every $i=0,1,\ldots,D$, the distance matrix $\A_i$ has entries $(\A_i)_{uv}=1$ if $\dist(u,v)=i$, and $(\A_i)_{uv}=0$ otherwise. In particular, $\A_0=\I$ and $\A_1=\A$.
From the positive (column) eigenvector $\vecalpha$, $\|\vecalpha\|^2=n$, we consider the matrices $\J^*=\vecalpha\vecalpha^{\top}$, with entries $(\J^*)_{uv}=\alpha_u\alpha_v$ for any $u,v\in V$, and $\A^*_i=\A_i\circ \J^*$, which can viewed as a ``weighted" version of the distance matrix $\A_i$. In fact, the approach of giving weights (which are the entries of the Perron vector) to the vertices of a nonregular graph has been recently often used in the literature (see, for instance, \cite{fgy96,fg97,fg99,f99,lw11}).

\subsection {Predistance polynomials}
In the vector space of symmetric matrices $\Re^{n\times  n}$, define the scalar product
$$
\langle \M,\N\rangle=\frac{1}{n}\tr (\M\N)=\frac{1}{n}\sum_{i,j}(\M\circ \N)_{ij}
$$
in such a way that, within the adjacency algebra of a graph $\G$ with adjacency matrix $\A$ and spectrum $\sp \G=\{\lambda_0^{m(\lambda_0)},\lambda_1^{m(\lambda_1)},\ldots,\lambda_d^{m(\lambda_d)}\}$, we can consider the following scalar product in $\Re^n[x]/(Z)$ (where $Z$ is the minimum polynomial of $\A$):
\begin{equation}\label{scalar-prod-glob}
\langle p,q\rangle_{\G}=\frac{1}{n}\tr (p(\A)q(\A))=\langle p(\A),q(\A)\rangle= \frac{1}{n}\sum_{i=0}^d m(\lambda_i)p(\lambda_i)q(\lambda_i) \qquad p,q\in \Re_{d}[x].
\end{equation}

Then, the {\em predistance polynomials}
$p_0,p_1,\ldots,p_d$, introduced in \cite{fg97} and extensively used in subsequent papers (they were first called the ``proper polynomials"; their present name was first proposed in \cite{f02}), are the orthogonal polynomials with respect to the product $\langle\cdot,\cdot\rangle_{\G}$, normalized in such a way that
$\|p_i\|_{\G}^2=p_i(\lambda_0)$. (This makes sense since, as it is well-known from the theory of orthogonal polynomials, $p_i(\lambda_0)>0$ for every $i=0,1,\ldots,d$, see, for instance, \cite{cffg09, vd08}.)
As every sequence of orthogonal polynomials, the predistance polynomials satisfy a three-term recurrence with coefficients being the {\em preintersection numbers}, a concept studied in \cite{ddfg10}.

There is a local version of the predistance polynomials which goes as follows. Given a vertex $u$ with $d_u+1$ different local eigenvalues, we can  consider the $u$-local scalar product
\begin{equation}\label{scalar-prod-loc}
\langle p,q \rangle_{u} = (p(\A)q(\A))_{uu}=\langle p(\A)\vece_{u},
q(\A)\vece_{u}\rangle=\sum_{i=0}^d m_{u}(\lambda_i)p(\lambda_i)q(\lambda_i),\quad p,q\in \Re_{d_{u}}[x],
\end{equation}
where we have used that $\vece_{u}=\sum_{i=0}^d \E_i\vece_{u}$ and $p(\A)\E_i=p(\lambda_i)\E_i$.
(Notice that the sum has at most $d_{u}$ nonzero terms.) Then, the {\em $u$-local predistance polynomials} $\{p_i^u\}_{0\le i\le d_{u}}$ are the orthogonal sequence with respect to the product $\langle\cdot,\cdot\rangle_{u}$, but now they are normalized in such a way that
$\|p^u_i\|_{u}^2=\alpha_u^2p^u_i(\lambda_0)$.
Notice that, since $\sum_{u\in V}m_u(\lambda_i)= m(\lambda_i)$, the above global scalar product (\ref{scalar-prod-glob}) is the mean over $V$ of the local products in (\ref{scalar-prod-loc}):
\begin{equation}
\label{mean-prod}
\langle p,q\rangle_{\G}=\frac{1}{n}\sum_{u\in V}\langle p,q\rangle_{u}.
\end{equation}

\subsection{Hoffman polynomials}
The {\it sum predistance polynomials} $q_j$, for $0\le j\le d$, are defined as $q_j=p_0+p_1+\cdots+p_j$, so that the {\it Hoffman  polynomial} $H$ \cite{hof63}, characterized by $H(\lambda_i)=n\delta_{0i}$ for $0\le i\le d$, turns out to be $H=q_d$ \cite{cffg09}, and satisfies
$H(\A)=\J$ if and only if $\G$ is regular.

Recently, some general (local or global) versions of the Hoffman polynomial has been considered (see, for instance,  \cite{fgy96,fg97,lw11}). Thus, with $\vecnu$ being any Perron vector, the {\em preHoffman polynomial} is defined as
$$
H= p_0+p_1+\cdots +p_d=\frac{\|\vecnu\|^2}{\pi_0}\prod_{i=0}^d (x-\lambda_i),
$$
where $\pi_0=\prod_{i=0}^d (\lambda_0-\lambda_i)$, and satisfies
$$
(H(\A))_{uv}=\nu_u\nu_v\qquad (u,v\in V).
$$
Similarly, if $\ev_u\G=\{\mu_0(=\lambda_0),\mu_1,\ldots,\mu_{d_{u}}\}$, the {\em $u$-local preHoffman polynomial} is
$$
H^u=p^u_0+p^u_1+\cdots+p^u_{d_{u}}=\frac{\|\vecnu\|^2}{\pi_0}\prod_{i=0}^{d_{u}} (x-\mu_i),
$$
where $\pi_0=\prod_{i=1}^{d_{u}} (\lambda_0-\mu_i)$.
Note that $H^u(\A)\vece_u=H(\A)\vece_u$ since $\vece_u $ has nonzero projections only on $\Ker (\A-\mu_i\I)$, $i=0,1,\ldots,d_u$.

%%%%%%%%%%%%%%%%%%%%%%%%%%%%%%%%%%%%%%%%%%%%%%%%
\subsection{Distance-regularity around a vertex}
\label{sec: d-r conjunts}
%%%%%%%%%%%%%%%%%%%%%%%%%%%%%%%%%%%%%%%%%%%%%%%%

Let $\G=(V,E)$ be a graph with  adjacency matrix $\A$, maximum eigenvalue $\lambda_0$ and Perron vector $\vecalpha$, $\|\vecalpha\|^2=n$.
Consider the map $\vecrho: V \longrightarrow \Re^n$ defined by
$\vecrho(u)= \vecrho_u=\alpha_u\e_u$, where $\e_u$ is the coordinate vector.
Note that, since $\|\vecrho_u\|=\alpha_u$, we can see $\vecrho$ as a function which assigns weights to the vertices of $\G$. In doing so we ``regularize" the graph, in
the sense that the {\it average weighted degree} of each vertex $u\in V$ becomes a
constant:
\begin{equation}
\label{regularize}
\delta^*_u=\frac{1}{\alpha_u}\sum_{v\in\G(u)}\alpha_v=\lambda_0.
\end{equation}

A graph $\G$ is said to be {\em pseudo-distance-regular around a vertex $u\in V$} with eccentricity $\ecc(u)=\excu$ (or {\em $u$-local pseudo-distance-regular\/}) if the numbers, defined for any vertex $v\in \G_i(u)$,
\begin{equation}\label{loc-intersec-num}
c_i^*(v)=\frac{1}{\alpha_v}\sum_{w\in\G_{i-1}(u)}\alpha_w,\quad
a_i^*(v)=\frac{1}{\alpha_v}\sum_{w\in\G_{i}(u)}\alpha_w,\quad
b_i^*(v)=\frac{1}{\alpha_v}\sum_{w\in\G_{i+1}(u)}\alpha_w,\quad
\end{equation}
depend only on the value of $i=0,1,\ldots,\excu$. In this case, we denote them by $c_i^*$, $a_i^*$, and $b_i^*$, respectively, and they are referred
to as the {\em $u$-local pseudo-intersection numbers} of $\G$.
%By (\ref{regularize}), notice that
%\begin{equation}
%\label{a++bk+ck}
%a^*_i+b^*_i+c^*_i=\lambda_0\qquad (0\le i\le \excu).
%\end{equation}
In particular, when $\G$ is regular, $\vecalpha=\j$ and the $u$-local pseudo-distance-regularity coincides with the distance-regularity around $u$ (see Brouwer, Cohen and Neumaier \cite{bcn89}).

\section{The {\em SPET\/} for nonregular graphs}
\label{sec-local}
The following result is stated
%In \cite{fg97} it was shown that the predistance polynomial $p^u_{d_{u}}$ (there it was normalized in a different way) satisfies
%$p^u_{d_{u}}(\lambda_0)\ge \|\vecrho_{\G_{d_{u}}(u)}\|^2$
%and equality holds if and only if $p^u_{d_{u}}(\A)\vece_{u}=\alpha_u\vecrho_{\G_{d_{u}}(u)}$.
%A more general result, stated
in terms of the {\em sum $u$-local predistance polynomials} $q_{j}^u=p_0^u+p_1^u+\cdots+p_{j}^u$.

\begin{proposicio}[\cite{fg99}]
\label{prop caract c-c-r sum pol}
Let $u$ be a vertex  of a graph, with $u$-local predistance polynomials $\{p_i^u\}_{0\le i\le d_{u}}$ and let $q_j^u=\sum_{i=0}^j p_i$. Then, for any polynomial $r\in \Re_j[x]$, $0\le j\le d_{u}$, we have
\begin{equation}
\label{c-c-r1}
\frac{r(\lambda_0)}{\|r\|_{u}}\leq \frac{1}{\alpha_u}\|\vecrho_{\subN_j(u)}\|,
\end{equation}
and equality holds if and only if $u$ is extremal and
\begin{equation}
\label{c-c-r2}
\frac{1}{\|r\|_{u}}r(\A)\vece_{u}=\vece_{\subN_j(u)}.
\end{equation}
In this case, $r=\eta q_j^u$ for any $\eta\in\Re$, whence $(\ref{c-c-r1})$ and $(\ref{c-c-r2})$ become, respectively,
\begin{equation}\label{c-c-r3}
q_j^u(\lambda_0)=\|\vecrho_{\subN_j(u)}\|^2,
\end{equation}
and
\begin{equation}\label{c-c-r4}
q_j^u(\A)\vece_{u}=\alpha_u \vecrho_{\subN_j(u)}\quad \iff \quad (q_j^u(\A))_{uv}=
\left\{
\begin{array}{ll}
\alpha_u\alpha_v & \mbox{if $\dist(u,v)\le j$,}\\
0 & \mbox{otherwise.}
\end{array}
\right.
\end{equation}
\end{proposicio}

As a consequence, the following characterizations of pseudo-distance-regularity around a vertex were also proved by the same authors, where $(b)$ can be seen as a {\em Local Spectral Excess Theorem}.

\begin{theo}[\cite{fg99}]
\label{teo caract c-c-r}
Let $\G=(V,E)$ be a connected graph and let $u$ be a vertex
 with eccentricity $\excu$. Then, $\G$ is pseudo-distance-regular around $u$ if and only if any of the following conditions holds:
\begin{itemize}
\item [$(a)$]
There exist polynomials $\{r_i\}_{0\leq i\leq \excu}$ with
$\dgr r_i=i$ such that
$$
r_i(\A)\vece_{u}=\alpha_u\vecrho_{\G_i(u)}.
$$
If this is the case, $u$ is extremal, $\excu=d_u$, and $r_i=p_i^u$ for $i=0,1,\ldots,d_u$.
\item [$(b)$]
The $u$-local predistance polynomial with highest
degree satisfies
$$
p_{d_u}^u(\lambda_0)=\|\vecrho_{\G_{d_u}(u)}\|^2.
$$
\end{itemize}
\end{theo}

Let $\delta^*_D=\|\A^*_D\|^2$ be the arithmetic mean of the numbers $\alpha_u^2\sum_{v\in\G_D(u)}\alpha_v^2$ for $u\in V$, which is called the {\em average weighted excess}, and let $p_{\ge D}(\lambda_0)=(q_d-q_{D-1})(\lambda_0)=n-q_{D-1}(\lambda_0)$ be the so-called {\em generalized spectral excess}. Then, Lee and Weng \cite{lw11} proved the following version of the spectral excess theorem for nonregular graphs.

\begin{theo}[\cite{lw11}]
\label{lee-weng-sp-ex-theo}
Let $\G$ be a connected graph with $n$ vertices, diameter $D$, weight distance matrices $\I^*,\A^*,\ldots,\A_D^*$, distinct eigenvalues $\lambda_0>\lambda_1>\cdots>\lambda_d$, predistance polynomials $p_0,p_1,\ldots,p_d$, and positive eigenvalue $\vecalpha$. Then,
\begin{equation}\label{lw-eq1}
\delta^*_D\le p_{\ge D}(\lambda_0),
\end{equation}
and equality is attained if and only if
\begin{equation}\label{lw-eq2}
\A^*_D=p_{\ge D}(\A).
\end{equation}
\end{theo}

The proof of this theorem is basically the same as the corresponding one in \cite{fgg10} (the only difference is the use of ``weights" on the vertices).
In \cite{lw11} Lee and Weng posed the problem of characterizing the graphs
for which the equality in (\ref{lw-eq1}) is attained, and gave an answer for two kinds of regular graphs. Namely,
when $D=d$, in which case the graphs are distance-regular; and when $D=2$, where the graphs turn to be distance-polynomial \cite{w82}. In fact, both cases had been already covered by the standard spectral excess theorem and some of its
generalizations \cite{ddfgg11}.
%%In fact, one could guess that equality in (\ref{lw-eq1}) cannot
%%occur unless the graph is regular.
%However, they did managed to prove the following theorem.
%\begin{theo}[\cite{lw11}]
%\label{theo-lee-weng}
%A connected graph $\G$ with $d+1$ distinct eigenvalues and odd-girth $2d+1$ is distance-regular.
%\end{theo}
%In fact, the same consequence had been previously obtained by Van Dam and
%Haemers \cite{vdh11} assuming the regularity of $\G$.
%(although these authors also proved Theorem \ref{theo-lee-weng} for the case
%$d+1=3$ and claimed to have proofs for the cases $d+1\in \{4,5\}$).

Now, we use the local approach of Proposition \ref{prop caract c-c-r sum pol} to obtain an improvement of the bound in (\ref{lw-eq2}), which leads to new characterizations of  some distance-regularity properties. We begin with the following more general result, where we use the sums of the weighted distance matrices $\S_j^*=\sum_{i=0}^j \A_i^*$, $j=0,1,\ldots,D$.

\begin{theo}
\label{teo(basic)}
Let $\G$ be a graph and, for some integer $j\le \min\{d_u: u\in V\}$, let $H_{\le j}^*$ be the harmonic mean of the numbers $\frac{1}{\alpha_u^2}\|\vecrho_{N_j(u)}\|^2$ for $u\in V$. Then,
$$
q_j(\lambda_0)\le H_{\le j}^*,
$$
and equality holds if and only if
 %there exist some Perron vectors $\vecnu^u$ for every $u\in V$
 %such that $(q_k(\A))_{uv}=\nu_u^u\nu_v^u$ if $\dist(u,v)\le k$ and $(q_k(\A))_{uv}=0$ otherwise.
 $$
 q_j(\A)= \S_j^*.
 $$
\end{theo}

\prova
With $r=q_j$,  the
inequality (\ref{c-c-r1}) of Proposition \ref{prop caract c-c-r sum pol} can be written as
\begin{equation}
\label{first-global-inequal}
\|q_j\|_u^2\ge \frac{\alpha_u^2 q_j(\lambda_0)^2}{\|\vecrho_{N_j(u)}\|^2}.
\end{equation}
Thus, by taking the average over all vertices, we have
$$
\frac{q_j(\lambda_0)^2}{n}\sum_{u\in V}\frac{\alpha_u^2}{\|\vecrho_{N_j(u)}\|^2}
\le
\frac{1}{n}\sum_{u\in V}\|q_j\|_u^2
=
\|q_j\|_{\G}^2
=
q_j(\lambda_0),
$$
where we used  (\ref{mean-prod}). Consequently,
\begin{equation}\label{q_j<=H}
q_j(\lambda_0)\le \frac{n}{\sum_{u\in V}\frac{\alpha_u^2}{\|\vecrho_{N_j(u)}\|^2}}=H_{\le j}^*
\end{equation}
Besides, equality can only hold if and only if all the inequalities in (\ref{first-global-inequal}) are also equalities and, hence, $q_j=\eta_u q_j^u$ for every vertex $u$ and respective constant $\eta_u$. Consequently, from (\ref{c-c-r4}),
\begin{equation}\label{c-c-r4-H}
q_j(\A)\vece_{u}=\eta_u \alpha_u \vecrho_{\subN_j(u)}\quad \iff \quad (q_j(\A))_{uv}=
\left\{
\begin{array}{ll}
\eta_u\alpha_u\alpha_v & \mbox{if $\dist(u,v)\le j$,}\\
0 & \mbox{otherwise.}
\end{array}
\right.
\end{equation}
To show that the constant $\eta_u$ does not depend on the vertex $u$,  let $u,v$ be two adjacent vertices. Then, by (\ref{c-c-r4-H}),
$$
\eta_u\alpha_u\alpha_v=(q_j(\A))_{uv}=(q_j(\A))_{vu}=\eta_v\alpha_v\alpha_u,
$$
and, hence $\eta_u=\eta_v$. Therefore, since  $\G$ is connected, we get the claimed result.
\final

A graph $\G$ with diameter $D$ is called {\em $m$-partially distance-regular}, for some $0\le m\le D$, if their predistance polynomials satisfy $p_i(\A)=\A_i$ for every $i\le m$, see \cite{ddfgg11}. It is immediate to show that every $m$-partially distance-regular with $m\ge 2$ must be regular. By using basically the same proof in \cite[Prop. 2.5]{ddfgg11}, we have the following result.

\begin{proposicio}
\label{propo1-m-part-drg}
Let $\G$ be a graph with $d+1$ distinct eigenvalues, sum predistance polynomials $q_j$,
and sums of weighted distance matrices $\S_j^*$, $j=0,1,\ldots,D$. Then,  $\G$ is $m$-partially distance-regular if and only if
\begin{equation}\label{eq-propo1-m-part-drg}
   q_j(\A)=\S_j^*\qquad \mbox{for}\qquad j=m-1,m.
\end{equation}
\end{proposicio}
In fact, in our ``weighted case" the proof gives that $p_i(\A)=\A_i^*$ for every $i=0,1,\ldots,m$ and, hence, we conclude  that the graph must be regular since $p_0(\A)=\I=\I^*$.

As an immediate consequence of the above and Theorem \ref{teo(basic)} we have:

\begin{proposicio}
\label{propo2-m-part-drg}
Let $\G$ be a graph and, as above, let $H_{\le j}^*$ be the harmonic mean of the numbers $\frac{1}{\alpha_u^2}\|\vecrho_{N_j(u)}\|^2$ for $u\in V$. Then,
\begin{equation}\label{eq-propo1-m-part-drg}
(q_{m-1}+q_{m})(\lambda_0)\le H_{\le m-1}^*+ H_{\le m}^*,
\end{equation}
and equality holds if and only if
$\G$ is regular and $m$-partially distance-regular
\end{proposicio}

As another corollary of Theorem \ref{teo(basic)}, by taking $j=D-1$ and considering that $p_{\ge D}=H-q_{D-1}$, we get the following result, which is an improvement of Theorem \ref{lee-weng-sp-ex-theo}.

\begin{theo}
\label{main-theo}
Let $\G$ be a graph with $H_{\le D-1}^*$ being the harmonic mean of the numbers $\frac{1}{\alpha_u^2}\|\vecrho_{N_{D-1}(u)}\|^2$ for $u\in V$, and let $\delta_D^*$  be the arithmetic mean of the numbers $\alpha_u^2\|\vecrho_{\G_{D}(u)}\|^2$. Then,
\begin{equation}\label{f-eq2}
p_{\ge D}(\lambda_0)\stackrel{(i)}{\ge} n-H_{\le D-1}^*\stackrel{(ii)}{\ge} \delta_D^*.
\end{equation}
Moreover, equality in $(i)$ holds if and only if
\begin{equation}\label{lee-weng}
p_{\ge D}(\A)=\A^*_D,
\end{equation}
whereas equality in $(ii)$ occurs when the numbers $\|\vecrho_{\G_{D}(u)}\|^2$ are the same for all $u\in V$.
\end{theo}
\prova
Inequality $(i)$ and the corresponding case of equality are consequences of Theorem  \ref{teo(basic)} with $k=D-1$:
\begin{equation}
\label{(i)}
p_{\ge D}(\lambda_0) =  n-q_{D-1}(\lambda_0)\ge n-H^*_{\le D-1} =n-\frac{n}{\sum_{u\in V}\frac{\alpha_u^2}{\|\vecrho_{N_{D-1}(u)}\|^2}}.
  \end{equation}
In order to prove $(ii)$, we have that  $\|\vecrho_{N_{D-1}(u)}\|^2=\|\vecalpha\|^2-\|\vecrho_{\G_{D}(u)}\|^2= n-\|\vecrho_{\G_{D}(u)}\|^2$. Then,
\begin{eqnarray}
\delta_D^* & = & \frac{1}{n}\sum_{u\in V}\alpha_u^2\|\vecrho_{\G_{D}(u)}\|^2
=\frac{1}{n}\sum_{u\in V}\alpha_u^2(n-\|\vecrho_{N_{D-1}(u)}\|^2)\nonumber \\
 & = & n-\frac{1}{n}\sum_{u\in V}\alpha_u^2\|\vecrho_{N_{D-1}(u)}\|^2 \label{(ii)}.
 \end{eqnarray}
 Thus, from (\ref{(i)}) and (\ref{(ii)}) we see that $n-H^*_{\le D-1}\ge \delta_D^*$ if and only if
$$
 n^2\le
 \sum_{u\in V}\alpha_u^2\|\vecrho_{N_{D-1}(u)}\|^2
 \sum_{u\in V}\frac{\alpha_u^2}{\|\vecrho_{N_{D-1}(u)}\|^2}
 $$
 or, using that $\|\vecalpha\|^2=\sum_{u\in V}\alpha_u^2=n$,
\begin{eqnarray*}
\sum_{u\in V}\alpha_u^2 \sum_{v\in V}\alpha_v^2  & = &
\sum_{u\in V}\alpha_u^4+ \sum_{\{u,v\}\subset V} 2\alpha_u^2\alpha_v^2\\
 & \le & \sum_{u\in V}\alpha_u^2\|\vecrho_{N_{D-1}(u)}\|^2\sum_{v\in V}\frac{\alpha_v^2}{\|\vecrho_{N_{D-1}(v)}\|^2}\\
 & = &\sum_{u\in V}\alpha_u^4+\sum_{\{u,v\}\subset V}\alpha_u^2\alpha_v^2\left(\frac{\|\vecrho_{N_{D-1}(u)}\|^2}{\|\vecrho_{N_{D-1}(v)}\|^2}+
 \frac{\|\vecrho_{N_{D-1}(v)}\|^2}{\|\vecrho_{N_{D-1}(u)}\|^2}\right),
\end{eqnarray*}
which is always true because, for any positive real numbers $r,s$, we have $\frac{r}{s}+\frac{s}{r}\ge 2$ with equality if and only if $r=s$. Therefore, equality in $(ii)$ happens if and only if $\|\vecrho_{N_{D-1}(u)}\|^2=\|\vecrho_{N_{D-1}(v)}\|^2$ for all $u,v\in V$, and the result follows from the fact that $\|\vecrho_{\G_{D}(u)}\|^2=n-\|\vecrho_{N_{D-1}(u)}\|^2$ for every $u\in V$.
\final

Notice that the invariance of the numbers $\|\vecrho_{\G_{D}(u)}\|^2$ follows also from (\ref{lee-weng}) by multiplying both terms by $\vecalpha$, and the common value is $p_{\ge D}(\lambda_0)$.

To compare the above bounds for $p_{\ge D}(\lambda_0)$, let us consider
the complete bipartite graph $K_{2,3}$, with spectrum $\sp K_{2,3}=\{\sqrt{6},0,-\sqrt{6}\}$ and Perron vector  $$\vecalpha=(\sqrt{5}/2,\sqrt{5}/2,\sqrt{5}/\sqrt{6},\sqrt{5}/\sqrt{6},
\sqrt{5}/\sqrt{6})^{\top},
$$ whence we obtain
$$
p_{\ge D}(\lambda_0)=\frac{3}{2}>n-H_{D-1}^*=\frac{25}{17}>\delta^*_D=\frac{35}{24}.
$$

As a consequence of Theorem \ref{main-theo}, we can now given a partial answer to the above mentioned problem posed in \cite{lw11}.

\begin{theo}
Let $\G$ be a graph with arithmetic means $\delta_i^*$ of the numbers $\alpha_u^2\|\vecrho_{\G_i(u)}\|^2$, $i=D-1,D$, and polynomials $p_{\ge D}$ and $p_{D-1}$ satisfying.
$$
\delta^*_D= p_{\ge D}(\lambda_0)\quad \mbox{and}\quad \delta^*_{D-1}= p_{D-1}(\lambda_0).
$$
Then $\G$ is distance-polynomial.
\end{theo}

\prova
By Theorem \ref{lee-weng-sp-ex-theo}, we have that $p_{\ge D}(\A)=\A^*_D$.
Moreover, from the hypotheses
\begin{eqnarray*}
q_{D-1}(\A) &=& (H-p_{\ge D})(\A)=\J^*-\A^*_D=\S_{D-1}^*, \\
q_{D-2}(\A) &=& (H-p_{\ge D}-p_{D-1})(\A)=\S_{D-1}^*-\A^*_{D-1}=\S_{D-2}^*.
\end{eqnarray*}
Hence, by Proposition \ref{propo1-m-part-drg}, $\G$ is regular and $(D-1)$-partially distance-regular. Consequently, $p_{\ge D}(\A)=\A_D$ and $\G$ is distance-polynomial, as claimed.
\final

\vskip 1cm

\noindent {\bf Acknowledgments.}
Research supported by the Ministerio de Educaci\'on y
Ciencia (Spain) and the European Regional Development Fund under
project MTM2011-28800-C02-01, and by the Catalan Research Council
under project 2009SGR1387.

%%%%%%%%%%%%%%%%%%%%%%%%%%%%%%%%%%%%%%%%%%%%%%%%
%Bibliografia
%%%%%%%%%%%%%%%%%%%%%%%%%%%%%%%%%%%%%%%%%%%%%%%%

\end{document}